\documentclass[11pt]{amsart}
\usepackage{amsfonts, amsmath, amssymb, amsthm}
\usepackage{epic}
\usepackage{tikz}
\usepackage{graphicx,color}
\usepackage{changebar}

\newtheorem{lemma}{Lemma}
\newtheorem{theorem}{Theorem}

\newtheorem{proposition}{Proposition}

\newtheorem{definition}{Definition}
\newtheorem{conjecture}{Conjecture}
\newtheorem{remark}{Remark}
\newtheorem{fact}{Fact}

\newtheorem{Acknow}{Acknowledgement}

\newcommand{\nth}[1]{$#1 \hbox{-th }$}

\newcommand{\ord}{{\rm ord} }

\newcommand{\rk}{{\rm rank } }

\newcommand{\Z}{\mathbb{Z}}
\newcommand{\ZM}[1]{\Z /( #1 \cdot \Z)}

\newcommand{\rg}[1]{\mathbf{#1}}
\newcommand{\eu}[1]{\mathfrak{#1}}
\newcommand{\re}[1]{\mathbf{\mathfrak{#1}}}
\newcommand{\id}[1]{\mathcal{#1}}
\newcommand{\Gal}{\hbox{Gal}}

\newcommand{\Ker}{\hbox{ Ker }}

\newcommand{\prk}{p\hbox{-rk} }
\newcommand{\eprk}{ess. \ p\hbox{-rk}}
\newcommand{\zprk}{\Z_p\hbox{-rk}}
\newcommand{\zrk}{\Z\hbox{-rk} }

\newcommand{\rf}[1]{(\ref{#1})}
\newcommand{\ran}{\rangle}
\newcommand{\lan}{\langle}
\newcommand{\Rad}{\hbox{Rad}}

\newcommand{\pinf}{^{1/p^{\infty}}}

\newcommand{\apr}[1]{{\rm \raise1ex\hbox{\tiny p}\kern-.1667em\hbox{A}}^-(#1)}
\newcommand{\apro}{{\rm \raise1ex\hbox{\tiny p}\kern-.1667em\hbox{A}}}
\newcommand{\ap}[1]{{\rm \raise1ex\hbox{\tiny p}\kern-.1667em\hbox{A}}(#1)}

\newcommand{\Norm}{\mbox{\bf N}}
 
\newcommand{\lchooses}[2]{\left( \frac{#1}{#2 } \right)}

\newcommand{\F}{\mathbb{F}}
\newcommand{\KP}{\mathbf{\Phi}}

\newcommand{\A}{\mathbb{A}}

\newcommand{\B}{\mathbb{B}}
\newcommand{\K}{\mathbb{K}}

\newcommand{\LK}{\mathbb{L}}
\newcommand{\KH}{\mathbb{H}}

\newcommand{\KL}{\mathbb{L}}

\newcommand{\M}{\mathbb{M}}
\newcommand{\Q}{\mathbb{Q}}

\newcommand{\C}{\mathbb{C}}
\newcommand{\N}{\mathbb{N}}

\def\ra{\rightarrow}

\def\wh{\widehat}

\begin{document}
{\obeylines \small
\vspace*{-1.0cm}
\hspace*{3.0cm}Par les gosses battus, par l'ivrogne qui rentre
\hspace*{3.0cm}Par l'\^{a}ne qui re\c{c}oit des coups de pied au ventre
\hspace*{3.0cm}Et par l'humiliation de l'innocent ch\^{a}ti\'{e}
\hspace*{3.0cm}Par la vierge vendue qu'on a d\'{e}shabill\'{e}e
\hspace*{3.0cm}Par le fils dont la m\`{e}re a \'{e}t\'{e} insult\'{e}e
\vspace*{0.2cm} 
\hspace*{4.0cm} Je vous salue, Marie\footnote{Francis Jammes: {\em Pri\`ere}. Music by Georges Brassens}
\vspace*{0.4cm}
\hspace*{4.5cm} {\it To Theres and Seraina, to the memory of Marta }
\vspace*{0.5cm}
\smallskip}
\title[Leopoldt's Conjecture] {On CM $\Z_p$-extensions and the Leopoldt conjecture for CM fields} \author{Preda
  Mih\u{a}ilescu} \address[P. Mih\u{a}ilescu]{Mathematisches Institut
  der Universit\"at G\"ottingen}
\email[P. Mih\u{a}ilescu]{preda@uni-math.gwdg.de}
\date{Version 2.0 \today}
\vspace{0.5cm}
\begin{abstract}
  We show that Leopoldt's conjecture holds in CM fields. For the proof
  we construct a $CM \Z_p$-extension of some CM field in which the
  Leopoldt conjecture is supposed to fail, and using the classes of
  primes which are completely split in this extension, we derive an
  contradiction. The method of proof can be described as a {\em
    stability check of $\Lambda$-modules under deformations in Thaine
    shifts}.
 \end{abstract}
\maketitle
\tableofcontents
\section{Introduction}
Leopoldt suggested in his seminal paper \cite{Le} that the $p$-adic
regulator of abelian extensions of $\Q$ never vanishes. This fact could
be proved by Brumer \cite{Br} in $1967$, using a plan of Ax \cite{Ax},
as soon as Baker had proved his archimedean version of the
approximation Theorem for linear forms in logarithms \cite{Ba}: it
remained to adapt Baker's proof to the $p$-adic topology. The fact
that Leopoldt's observation should hold in arbitrary number fields
became soon a largely accepted conjecture. It should however be noted,
that Leopoldt's definition of the $p$-adic regulator as a formal 
adaptation of the global regulator, in which logarithms are taken $p$-adically,
only makes sense in totally real or CM fields. For non-CM fields, there is 
no simple and canonical definition known: therefore, the Leopoldt conjecture
for this case should be related to the equality of $\Z$- and $\Z_p$ ranks of
the units, as defined below.
The purpose of this paper is to prove:

\begin{theorem}
\label{main}
For odd primes $p$, the Leopoldt Conjecture holds in arbitrary CM
extensions $\K/\Q$.
\end{theorem}

\subsection{Notations and fundamental facts}
Let $p$ be an odd rational prime. In this paper we denote
number fields with black board bold letters -- $\B, \F, \KL, \Q$, etc.
Bold characters will be used in general for local fields or for
denoting generic fields in some general remarks. We consider
CM number fields $\K$; their cyclotomic $\Z_p$-extension
will be denoted by $\K_{\infty}$, while $\KL$ will denote
some other $\Z_p$-extension of the same field. 
The intermediate fields will be $\K_n$, resp. $\KL_n$.

If $\KL/\K$ is an arbitrary $\Z_p$-extension, with group $\Gamma \cong
\Z_p$, generated by $\tau \in \Gamma$ as a topological generator, the
Iwasawa algebra is $\Lambda = \Z_p[[ T ]], T = \tau - 1$. The
intermediate fields are $\KL_n \subset \KL$ and, if $\KL =
\K_{\infty}$ is the cyclotomic $\Z_p$-extension, then we always assume
that $\K_n = \KL_n$ contains the \nth{p^n} roots of unity, but
not the \nth{p^{n+1}} ones; this can be achieved by an adequate
numeration, at least for some sufficiently large $n$, as shall be shown
in more detail in the next chapter. The $\Z_p$-extension of $\Q$
is denoted by $\B_{\infty}$ and we count the intermediate fields according to
$\B_n = \Q[ \zeta_{p^n} ] \cap \B_{\infty}$, with $\zeta_{p^n}$ some primitive
\nth{p^n} root of unity.
 
We write $\tau_n = (T+1)^{p^{n-k}}$ for the power of $\tau$ that
generates the fixing group of $\K_n$, and $\omega_n = \tau_n - 1$;
for $j > 0$, we let $\nu_{n+j,n} = \omega_{n+j}/\omega_n$. The $p$-parts of the
ideal class groups of $\K_n, \KL_n$ are $A_n = A(\K_n), A(\KL_n)$.
The \textit{injective} limit is denoted traditionally by $A = \varinjlim_n A_n$,
but we shall not use this limit in the present paper. 

Let $\KH(\K_n)$ be the maximal $p$-abelian unramified extensions of
$\K_n$, \ for all $n \in \N$ and let $X_n = X(\K_n) =
\Gal(\KH_n/\K_n)$. Furthermore, $\KH(\K) = \cap_{n = 1}^{\infty}
\KH(\K_n)$ and $X = \varprojlim_n X_n$, the projective limit being
taken with respect to restriction maps.  The module $X$ is a
Noetherian $\Lambda$-torsion modules on which complex conjugation
acts, inducing the decomposition $X = X^+ \oplus X^-$.  The Artin maps
$\varphi : A_n \ra X_n$ are isomorphisms of finite abelian $p$-groups;
whenever the abelian extension $\M/\K$ is clear in the context, we
write $\varphi(a)$ for the Artin Symbol $\lchooses{\M/\K}{a}$, where
$a \in \id{C}(\K)$ if $\M$ is unramified, or $a$ is an ideal of $\K$
otherwise.  Complex conjugation acts on class groups and Galois
groups, inducing direct sum decompositions in plus and minus parts:
\[ A_n = A_n^+ \oplus A_n^-, \quad X_n = X_n^+ \oplus X_n^-, \quad
\hbox{etc.} \] The idempotents which act on abelian groups, generating
plus and minus parts are $\frac{1 \pm \jmath}{2}$; since $2$ is a unit
in the ring $\Z_p$ acting on these groups, we also have $Y^+ =
(1+\jmath) Y$ and $Y^- = (1-\jmath) Y$, for $Y \in \{ A_n, X_n, X,
\ldots \}$.  Throughout the paper, we shall use, by a slight abuse of
notation and unless explicitly specified otherwise, the additive
writing for the group ring actions. This is preferable for
typographical reasons.

Our {\em working base field} $\K$ will be endowed with a series of
properties which shall be discussed in the next chapter. In
particular, the norms $\Norm_{\K_m/\K_n} : A_m \ra A_n$ are surjective
for all $m > n > 0$. The sequence $(A_n)_{n \in \N}$ is projective
with respect to the norm maps and we denote their projective limit by
$\apro = \varprojlim_n A_n$. The Artin map induces an isomorphism of
compact $\Lambda$-modules $\varphi :\ \apro \ra X$. The elements of
$\apro$ are norm coherent sequences $a = (a_n)_{n \in \N} \in \apro$
with $a_n \in A_n$ for $n \geq \kappa$; we let $a_0 = a_1 = \ldots =
a_{\kappa}$. It is customary to identify $X$ with $\apro$ via the
Artin map.

Let $M$ be a Noetherian $\Lambda$-torsion module. It is associated to
an \textit{elementary} Noetherian $\Lambda$-torsion module $\id{E}(M)
\sim M$ defined by:
\begin{eqnarray*}
\begin{array}{c c c c c c c }
\id{E}(M) & = & \id{E}_{\lambda}(M) & \oplus & \id{E}_{\mu}(M), & \quad & \hbox{with}   \\
\id{E}_{\mu}(M) & = & \oplus_{i=1}^r \Lambda/(p^{e_i} \Lambda), & \quad & \id{E}_{\lambda}(M) & = & \oplus_{j=1}^{r'} 
\Lambda/(f_j^{e'_j} \Lambda),  
\end{array}
\end{eqnarray*}
where all $e_i, e'_j > 0$ and $f_j \in \Z_p[ T ]$ are irreducible
distinguished polynomials.  The pseudoisomorphism $M \sim \id{E}(M)$
is given by an exact sequence
\begin{eqnarray}
\label{psis}
  1  \ra  K_1  \ra  M  \ra  \id{E}(M)  \ra  K_2  \ra  1 ,  
\end{eqnarray}
in which the kernel and cokernel $K_1, K_2$ are finite. We define
$\lambda$ and $\mu$-parts of $M$ as follows:
\begin{definition}
  Let $M$ be a Noetherian $\Lambda$-torsion module. The $\lambda$-part
  $\id{L}(M)$ is the maximal $\Lambda$-submodule of finite $p$-rank
  contained in $M$. The $\mu$-part $\id{M}(M)$ is the $\Z_p$-torsion
  submodule of $M$; it follows from the Weierstra{\ss} Preparation
  Theorem that there is some $m > 0$ such that $\id{M}(M) = M[ p^m
  ]$. The maximal finite $\Lambda$-submodule of $M$ is $\id{F}(M)$,
  its finite part. By definition, $\id{L}(M) \cap \id{M}(M) =
  \id{F}(M)$.

  Let the module $\id{D}(M) = \id{L}(M) + \id{M}(M)$ be the {\em
    decomposed submodule} of $M$.  Then for all $x \in \id{D}(M)$
  there are $x_{\lambda} \in \id{L}(M), x_{\mu} \in \id{M}(M)$ such
  that $x = x_{\lambda} + x_{\mu}$, the decomposition being unique iff
  $\id{F}(M) = 0$.  The pseudoisomorphism $M \sim \id{E}(M)$ implies
  that $\id{D}(M) \sim \id{E}(M)$ too, and thus $[ M : \id{D}(M) ] <
  \infty $.

  If $x \in M \setminus \id{D}(M)$, the $L$- and the $D$-orders of $x$
  are, respectively
\begin{eqnarray}
\label{dords}
\ell(x) & = & \min\{  j > 0 \ : \ p^{j} x \in \id{L}(M) \}, \quad \hbox{and} \\
\delta(x) & = & \min\{  k > 0 \ : \ p^{k} x \in \id{D}(M) \} \leq \ell(x). \nonumber
\end{eqnarray}
\end{definition}


If $\K$ is a number field, we denote its units by $E(\K) =
\id{O}^{\times}(\K)$. Dirichlet's unit theorem states that, up to
torsion made up by the roots of unity $W(\K) \subset \K^{\times}$, the
units $E = \id{O}(\K)^{\times}$ are a free $\Z$ - module of $\Z$ -
rank $r_1 + r_2 - 1$. As usual, $r_1$ and $r_2$ are the numbers of
real, resp. pairs of complex conjugate embeddings $\K \hookrightarrow
\C$. We consider the set $P = \{ \wp \subset \id{O}(\K) : (p) \subset
\wp \}$ of distinct prime ideals above $p$ and let
\[ \re{K}_p = \re{K}_p(\K) = \prod_{\wp \in P } \K_{\wp} = \K \otimes_{\Q} \Q_p
\] 
be the product of all completions of $\K$ at primes above $p$.  Let
$\iota : \K \hookrightarrow \re{K}_p$ be the diagonal embedding. We
write $\iota_{\wp}(x)$ for the projection of $\iota(x)$ in the
completion at $\wp \in P$. If $y \in \re{K}_p$, then $\iota_{\wp}(y)$
is simply the component of $y$ in $\K_{\wp}$. If $U \subset
\re{K}_p^{\times}$ is the group of units, thus the product of local
units at the same completions, then $E$ embeds diagonally via $\iota :
E \hookrightarrow U$.

Let $\overline{E} = \overline{\iota(E)} = \bigcap_{ n > 0} \iota(E)
\cdot U^{p^n} \subset U$ be the $p$-adic closure of $\iota(E)$; this
is a $\Z_p$ - module with $\zprk(\overline{E}) \leq \zrk(E) = r_1 +
r_2 - 1$. The difference
\[ \id{D}_L(\K) = \zrk(E) - \zprk(\overline{E}) \] is called the
\textit{Leopoldt defect}. The defect is positive if, in the idelic
embedding, units which are $\Z$-independent, are related $p$-adically.
In CM fields, this is equivalent to the vanishing of 
the $p$ - adic regulator of $\K$. 

The connection of Leopoldt's conjecture to class field theory was
noted by Iwasawa in \cite{Iw}. He shows that if $\Omega(\K)
\supset \K_{\infty}$ is the maximal $p$-abelian $p$-ramified extension
of $\K$, then $\Gal(\Omega(\K)/\K) \sim \Z_p^n$, where $n =
r_2+1+\id{D}_L(\K)$; the proof of this fact is in any text book on
cyclotomy and Iwasawa theory. For CM extensions $\K$, the
contraposition of the conjecture herewith reduces to the statement
that $\K^+$ has more $\Z_p$-extensions than just the cyclotomic
one. It is this assumption which we shall use and lead to a
contradiction. If Leopoldt fails thus for $\K^+$ and if $\KL^+/\K^+$
is a further $\Z_p$-extension, it will be totally real and $\KL =
\KL^+ \cdot \K$ will be a CM extension. For such fields, an equivalent
formulation of the Leopoldt conjecture is thus:
\begin{conjecture}[ Leopoldt ]
\label{le}
The only $\Z_p$-extension of a totally real field is the cyclotomic
one, and a CM number field has no CM $\Z_p$-extensions except for the
cyclotomic one.
\end{conjecture}
We shall use this particular formulation of the conjecture in this
paper.
\subsection{Historical notes}
Since 1967 various attempts have been undertaken for extending the
results of \cite{Br} to non abelian extensions, using class field
theory, Diophantine approximation or both. The following very succinct
list is intended to give an overview of various approaches, rather
than being an extensive list of results on Leopoldt's conjecture. In
\cite{Gr}, Greenberg notes a relation between the Leopoldt Conjecture
and a special case of the Greenberg Conjecture: he shows that
Leopoldt's Conjecture implies that the $T$-part $\apro(T)$ is finite
for totally real fields, i.e. the Greenberg Conjecture holds for the
$T$ - part.

Emsalem, Kisilevsky and Wales \cite{EKW} use group representations and
Baker theory for proving the Conjecture for some small non abelian
groups; this direction of research has been continued in some further
papers by Emsalem or Emsalem and coauthors. Jaulent proves in
\cite{Ja2} the Conjecture for some fields of small discriminants,
using the \textit{phantom} field $\KP$ which we define in the
Appendix. Currently the strongest result based on Diophantine
approximation was achieved by Waldschmidt \cite{Wal}, who proved that
if $r$ is the $\Z$ - rank of the units in the field $\K$, then the
Leopoldt defect satisfies $\id{D}_L(\K) \leq r/2$.

\subsection{Plan of the proof}
We start with several general observations.
\begin{remark}
\label{shift}
\begin{itemize}
\item[ A. ] If $\K_{start}$ is a field for which $\id{D}_L(\K_{start})
  > 0$, then it is known that the same holds for arbitrary finite
  algebraic extensions $\K/\K_{start}$; this is noted, for instance,
  by Laurent in the introduction to \cite{Lau}. It follows from the
  fact that the linear relations between $\Z$-generators of the units
  of $E(\K_1)$, which arise upon $p$-adic completion, will be
  preserved under the embedding into the units $E(\K)$.

\item[ B. ] If $\K$ is a CM extension containing the \nth{p^m} roots
  of unity and $\KL^+/\K^+$ is a $p$-ramified cyclic extension of
  degree $p^m$, then there is a class $a \in A^-(\K)$ such that $\KL^+
  \cdot \K = \K[ a^{1/p^m} ]$ in the sense that for each $\eu{A} \in
  a$ there is an $\alpha$ generating the principal ideal $\eu{A}^{(1-\jmath) p^m}$,
  and such that $\KL = \K[ \alpha^{1/p^m} ]$. We shall use this notation for explicit Kummer
  extensions throughout the paper.

\item[ C. ] The base field $\K$ can be modified within the same
$\Z_p$-extension, by replacing $\K$ with $\K_n$, say. As a
consequence, the Iwasawa algebra may become $\Lambda' = \Z_p[[ \omega_n ]]$.
\end{itemize}
\end{remark}

The proof is inspired by a construction of Iwasawa for showing that
there exist $\Lambda$-extensions with $\mu > 0$. We assume that
$\K_{start}$ is some CM extension for which the Leopoldt defect does
not vanish, and create first an auxiliary extension $\K \supset
\K_{start}$ which also has positive Leopoldt defect, together with
some additional properties. This construction of a working base field
will be the topic of the first part of the second chapter. The working
base field $\K$ comes together with a norm coherent sequence
$\tilde{a} = (\tilde{a}_m)_{m \in \N} \in \apr{\K_{\infty}}$.  For
some sufficiently large $N$, we construct a CM $\Z_p$-extension $\KL/\K$
which intersects the cyclotomic one in $\K_{N} = \KL \cap
\K_{\infty}$ and has the property that a chosen prime $\eu{q} \in
\tilde{a}_{N}$ is totally split in $\KL/\K_{N}$. Let the
intermediate fields of $\KL$ be $\KL_n$, with $\KL_{N} = \K_{N}$, as
previously noted. We choose a norm coherent sequence of primes
$\eu{q}_n \subset \id{O}(\KL_n)$ above $\eu{q}$ and let $a_n = [
\eu{q}^{1-\jmath}_n ]$ be their classes.

Using the same prime $\eu{q}$ we construct a {\em split Thaine shift}:
this is the compositum of $\KL$ with the subfield $\F \subset \Q[
\zeta_q ]$ of degree $p$ in the \nth{q} cyclotomic field, where $q$ is the
rational prime below $\eu{q}$. We let $\KL_n' = \KL_n \cdot \F$ and
$\eu{Q}_n \subset \KL'_n$ be the ramified primes above $\eu{q}_n$.
They define the norm coherent sequence of classes $b_n := [
\eu{Q}_n^{1-\jmath} ]$ which form the {\em Thaine lift} of $a_n$; let
$a = (a_n)_{n \in \N} \in \apr{\KL}$ and $b := (b_n)_{n \in \N} \in
\apr{{\KL'}}$. Due to ramfication, the lift is in this case unique and
verifies $a = p b$. A cohomological investigation, using the Hasse Norm
Principle, shows that $b$ must be indecomposed. With this
construction in place, it is then a matter of a straightforward
algebraic computations, to obtain a contradiction with the premises
about the choice of $\tilde{a}$ and $\eu{q}$, which yields the proof
of the Conjecture. This proof is provided in the third chapter. 

\section{Auxiliary constructions}
In this chapter we construct the base working field and the Thaine shift 
mentioned in the introduction. Let $\K_{start}$ be a CM number field 
in which the Leopoldt conjecture fails for the odd prime $p$. 
In view of point A. in Remark \ref{shift}, we shall start by constructing
a sequence of extensions of $\K_{start}$, eventually obtaining a CM extension 
$\K \supset \K_{start}$ endowed with a series of properties that are required for
our proof, and in which the Leopoldt conjecture is thus false, too. We start
with the construction of an auxiliary field, which is essentially independent
of $\K_{start}$.

We adopt throughout this paper the following {\em counting
  convention}. Let $\M$ be an arbitrary number field and $\B_{\kappa}
= \M \cap \B_{\infty}$.  We then count the intermediate fields of the
cyclotomic $\Z_p$-extension of $\M$ as follows: $\M = \M_0 = \M_1 =
\ldots \M_{\kappa} \subsetneq \M_{\kappa+1}$.  This way, for $m >
\kappa$ we always have $\M_m \cap \B_{\infty} = \B_m$. All groups
associated to the intermediate fields are counted accordingly, thus
$A_0(\M) = \ldots = A_{\kappa}(\M)$.  In addition, if $\M'/\M$ is some
other $\Z_p$-extension which intersects the cyclotomic nontrivially,
say $\M' \cap \M_{\infty} = \M_n$, then the counting of intermediate
fields of $\M'$ is induced from the one for $\M_{\infty}$. Thus
$\M'_{n+1} \subset \M'$ will be the subfield of degree $p$ over $\M_n
= \M'_n$.

We start with the investigation of the torsion of minus parts in CM
$\Z_p$-extensions.  The first two sections of this chapter provide
some general facts about CM fields and CM $\Z_p$-extensions
thereof. They are herewith dependent on the failure of the Leopoldt
conjecture for the base field -- but its choice is still arbitrary,
except for the fact that it should possibly be galois and contain the
\nth{p} roots of unity.

\subsection{On CM $\Z_p$-extensions of number fields and auxiliary
  fields}
In this section we gather several properties of CM $\Z_p$-extensions
which are the base for our approach; recall that the occurrence of CM
$\Z_p$-extensions different from the cyclotomic one, is
\textit{equivalent} to the failing of Leopoldt's conjecture for CM
fields $\K$. We let $\K$ be some galois CM number field for which the
Leopoldt conjecture fails and let $\K_{\infty}$ be its cyclotomic
$\Z_p$-extension. We let $\M$ be the compositum of all the
$\Z_p$-extensions\footnote{Oberve that $\M$ was used previously to denote
some generic CM extension, while it is from now on the explicit product of all
CM extension of a chosen field $\K$.} of $\K$, let $\M_0^+$ be the compositum of all the
$\Z_p$-extensions of $\K^+$ and $\M^+ = \K \cdot \M_0^+ \subset \M$.

We denote by $\apr{\K}(T^*) = \{ x \in \apr{\K} \ : \ (T^*)^M x = 0 \hbox{ for some $N > 0$} \}$
the $T^*$-part\footnote{ Here $* : \Lambda \ra \Lambda$ 
denotes the Iwasawa involution induced by $\tau \mapsto \chi(\tau) \tau^{-1}$, 
with $\chi$ the cyclotomic character} of $\apr{\K}$. The $T^*$-torsion is $\apr{\K}[ T^* ]$.
The radicals of $\M^+$ as a Kummer extension of $\K_{\infty}$ are
intimately related to the failure of Leopoldt's conjecture and to the
$T^*$-part of the class group, by the following \textit{folklore}
result, which holds in the cyclotomic $\Z_p$-extension of a field:
\begin{proposition}
\label{fantom}
Let $\K$ be a CM field which contains the \nth{p} roots of unity and
$\ap{\K} = \varprojlim_n A(\K_n)$ be defined with respect to the
cyclotomic $\Z_p$-extension. Then
\[ \zprk(\apro[ T^* ]) = \id{D}_L(\K) , \] and in particular Leopoldt's
conjecture fails for $\K$ iff $\apro[ T^* ] \neq 0$. Moreover
\begin{eqnarray}
\label{radfan}
 \M^+ \subseteq \K_{\infty}[ (\apro( T^* ))\pinf ].
\end{eqnarray}
In particular, for every CM $\Z_p$-extension $\KL/\K$ there is a class
$a \in \apro( T^*) \setminus T^* \apro(T^*)$ such that $\KL \cdot
\K_{\infty} = \K_{\infty}[ a\pinf ]$.
\end{proposition}
The proof of the proposition is sketched in the Appendix.

For the cyclotomic $\Z_p$-extension, it is known that
$\apr{\K_{\infty}}$ has no finite $p$-torsion submodule. In the case of
non-cyclotomic CM $\Z_p$-extensions, this fact is almost true, namely:
\begin{lemma}
\label{almin}
Let $\K$ be a CM extension containing the \nth{p} roots of unity and
$\KL/\K$ be a CM $\Z_p$-extension with $\KL \cap \K_{\infty} = \K_N, N
\geq 1$ and write $\KL_N := \K_N; \ [ \KL_{N+n} : \K_N ] = p^n$. If
$\mu_{p^N} \subset \KL$ but $\mu_{p^{N+1}} \not \subset \KL$, then the
finite torsion submodule $C^- := \id{F}(\apr{\KL}) \subset
\apr{\KL}$ is a cyclic group of order $p^N$ annihilated by $T^*$.
\end{lemma}
\begin{proof}
  Since the norms $N_{m,n} : A^-(\KL_m) \ra A^-(\KL_n)$ are surjective
  we have $| C^- | \geq | C_m^- | \geq | C_n^- |$. The group sizes $|
  C_n^- |$ build thus an increasing, bounded sequence of integers, so
  there is some $n_0$ such that $| C^-_n | = | C^- |$ for all $n >
  n_0$; in fact, there are isomorphisms $C_n^- \cong C^-$ for all such
  $n$. For such $n$, we have also $\nu_{n+1,n} (C_{n+1}^-) = p \cdot C_{n+1}^-$ for $n > n_0$.
  Since the lift $\iota_{n,n+1}(C^-_n) = \nu_{n+1,n} C^-_{n+1}$, it
  follows that $ \iota_{n,n+1}(C^-_n) = p C^-_{n+1}$. The isomorphism
  $C^-_{n} \cong C^-_{n+1}$ implies that the capitulation kernel $P_n
  := \Ker(\iota_{n,n+1} : C^-_n \ra C^-_{n+1})$ has maximal rank:
  $\prk(P_n) = \prk(C_n^-) =\prk(C^-)$. Moreover, $N_{n+1,n} \circ
  \iota_{n,,n+1} = \cdot p$ implies that $\exp(P_n) = p$ and thus $P_n
  = C^-_n[ p ]$ for all $n > n_0$.

  We claim that $\prk(P_n) = 1$, which implies that $C^-$ must be
  $\Z_p$-cyclic.  Let $c \in C_n^-[ p ]$ and $\eu{C} \in c$ be a prime
  with $\iota_{n,n+1}(\eu{C}) = (\delta)$ and $\eu{C}^p = (\gamma), \
  \gamma \in \KL_n^{\times}$ and $\delta \in \KL_{n+1}^{\times}$.
  Then there is a unit $\varepsilon \in E(\KL_{n+1})$ such that
  $\gamma = \varepsilon \cdot \delta^p$.  By dividing through the
  complex conjugates, we deduce that there is some root of unity $\xi
  \in \mu_{p^N}$ such that $\gamma/\overline{\gamma} = \xi^2 \cdot
  (\delta/\overline{\delta})^p$.  Letting $\gamma' = \overline{\xi}
  \cdot \gamma$, we find that
  \[ \frac{\gamma'}{\overline{\gamma}'} =
  \left(\frac{\delta}{\overline{\delta}}\right)^p . \] The extension
  $\KL_{n+1}/\KL_n$ is a cyclic Kummer extension of degree $p$, so let
  $\beta \in \KL_n$ be such that $\KL_{n+1} =
  \KL_n\left[ (\beta/\overline{\beta})^{1/p} \right]$.  We deduce by Kummer theory that
  \[ \frac{\gamma'}{\overline{\gamma}'} = w^p \cdot
  \left(\frac{\beta}{\overline{\beta}}\right)^a, \quad \hbox{ for some
    $w \in \KL_n^{\times}$ and $(a, p) = 1$} . \] In terms of point
  2. in Remark \ref{shift}, we have thus $\KL_{n+1} = \KL_n[ c^{1/p}
  ]$. This holds for arbitrary $c \in C_n^-[ p ]$, implying that
  $C_n^-[ p ]$ must be cyclic, as claimed.

  Let now $p^j = \exp(C^-) = | C^- |$. By an argument similar to the
  one used above, we deduce that for arbitrary $c \in C^-_n \setminus
  (C_n^-)^p$ we have $\KL_{n+j} = \KL_n[ c^{1/p^j} ]$.  In addition,
  $\KL/\K_N$ is abelian, so a straight forward computation implies
  that $p^j = | W(\KL) | = p^N$, so $j = N$ and $C^- \cong
  \ZM{p^N}$. 
  
  Finally, $\KL/\K$ is by definition a $\Z_p$-extension and thus $T
  \Gal(\KL/\K) = 0$. In particular, $T \Gal(\KL_{n+N}/\KL_n) = 0$ for
  all $n > n_0$, and Kummer duality implies that $T^* c = 0$, which
  completes the proof.
\end{proof} 

We note that there is some sequence $\beta = (\beta_n)_{n \in \N} \in
\apr{\K_{\infty}}(T^*)$, such that $\KL \cdot \K_{\infty} =
\K_{\infty}[ \beta\pinf ]$. One can verify that the torsion of
$\apr{\KL}$ is built from the restriction of classes in $\beta$ to the
intermediate fields of $\KL$, which also explains that the torsion is
annihilated by $T^*$.

We have shown in \cite{MiMu} that the Iwasawa constant $\mu$ vanishes
for cyclotomic $\Z_p$-extensions of CM fields. This result simplifies
the construction of our base field and shall be applied here. We
provide however in the Appendix a construction which does not depend
on the vanishing on $\mu$ and can thus be used in contexts in which
one wishes to avoid use of the proof in \cite{MiMu}.

We know thus that $\mu(\apr{\K_{\infty}} ) = 0$, so $\apr{\K_{\infty}}
= \id{L}(\apr{\K_{\infty}} )$.  Let $f(T)$ be the minimal annihilator
polynomial of $\apr{\K_{\infty}}$ and let $n_0$ be chosen such that
the following conditions hold:
\begin{itemize}
\item[ N1. ] For all $x = (x_n)_{n \in \N} \in \apr{\K_{\infty}}
  \setminus (p, T) \cdot\apr{\K_{\infty}}$ we have $\ord(a_{n_0}) >
  p^2$.
 \item[ N2. ] The polynomial $\omega_{n_0} \in (p, f(T)) \Lambda$.
\end{itemize}

The condition N1 can be fulfilled, since $M := \apr{\K_{\infty}} / ((p, T) \cdot\apr{\K_{\infty}})$ has finite rank
and for each $x$ with non trivial image in $M$, the orders $\ord(x_n)$
diverge, so there is a minimal $n'_0$ for which condition
N1. holds. For N2, notice that $(p, F(T)) \supset (p, T^{\lambda})$,
where $\lambda = \prk(\apr{\K_{\infty}})$ and $\omega_k \equiv
T^{p^{k-\kappa}} \bmod p$, so it suffices to choose $n_0 \geq \max
(n_0', \kappa + v_p(\lambda) + 1)$. We show in the Appendix, that
if $\mu > 0$ is assumed, it will suffice to replace condition N2 by 
a sharper condition, which also depends solely on $\K$. 

\begin{lemma}
\label{mubd} 
There is some $N > n_0$, such that for any CM $\Z_p$-extension $\M'/\K$ in with
$\M' \cap \K_{\infty} \supset \K_N$, we have $\mu^-(\M') = 0$.
\end{lemma}
\begin{proof}
  Note first that for $n > k$ all the primes above $p$ are totally
  ramified in $\K_n/\K$ and a fortiori in $\M'/\K$. Assume now that
  $\mu^-(\M') > 0$ and let $x \in \id{M}(\apr{\M'}) \setminus (p,T)
  \id{M}(\apr{\M'})$.  Then $x_{k+1} \neq 0$: otherwise, Iwasawa's
  Theorem VI \cite{Iw} would imply that there is some $y \in \apr{\M'}$
  such that $x = \nu_{k+1,1} \ y$ and thus $\ord(x) \cdot y \in
  \Ker(\nu_{k+1,1}) = 0$.  Therefore $y \in \id{M}(\apr{\M'})$, in
  contradiction with the choice of $x$. Since $x_{k+1} \neq 0$, and $x
  \in \id{M}(\M')$, it follows that $\prk(\Lambda x_n) \geq
  p^{n-{k+1}}$ for $n > k$. However, $\prk(\Lambda x_n) \leq
  \lambda(\K_{\infty}) < \infty$, so by choosing $N = n$ sufficiently
  large, we obtain a contradiction, which implies that $\mu^-(\M') =
  0$.
\end{proof}

We show how to construct an abelian extension of $\Q$, which will be an
auxiliary step in the construction of our working base field
\begin{lemma}
\label{paux}
Let $m, k > 0$ be integers. There is an imaginary quadratic extension
$\rg{K}$ which is linearly disjoint from $\K$ and does not contain the
cube roots of unity, and an abelian extension $\K_{aux}/\Q$ with
$\rg{K} \subset \K_{aux}$, which is unramified at the primes that
ramify in $\K$ and such that $\apr{\K_{aux, \infty}} \neq
0$. Moreover, there is a norm coherent sequence
\begin{eqnarray}
\label{palpha}
\alpha & = & (\alpha_n)_{n \in \N} \in \nonumber
\apr{\K_{aux, \infty}} \setminus (p, T) \cdot \apr{\K_{aux, \infty}}, \quad \hbox{such that}\\
\omega_m^k \alpha & \not \in & p  \apr{\K_{aux, \infty}}.
\end{eqnarray}
\end{lemma}
\begin{proof}
  Let $\rg{K} = \Q[ \sqrt{-D} ]; D \neq 3$ be an imaginary quadratic extension
  which is unramified at the rational primes which ramify in $\K$ and
  let $m'$ be such that $\deg(\nu_{m',1}) > \deg(\omega_m^{k+1})$ as
  polynomials in $T$. Note that $\rg{K}$ satisfies the premises of the Lemma.
  Let $\eu{r} \subset \rg{K}_{m'}$ be a totally
  split principal ideal above a rational prime $r \equiv 1 \bmod p$
  which does not ramify in $\K$, and let $\F' \subset \Q[ \zeta_r ]$
  be the subfield of degree $p$ in the \nth{r} cyclotomic extension.
  We assume in addition that $\eu{r}$ is inert in
  $\rg{K}_{\infty}/\rg{K}_{m'}$. Then $\rg{L} = \rg{K} \cdot \F'$ is
  the {\em inert Thaine shift} of $\rg{K}$ through the prime $\rg{r}$,
  and was studied in \S 2 of \cite{MiMu}. Let $R = \F_p[
  \Gal(\rg{K}_{m'} / \rg{K}) ]$ and $F' = \Gal(\F'/\Q)$. We apply the Proposition \ref{genHasse}
 below, in which we let $n = m'$ and $\KL = \rg{K}, \KL' = \rg{L}$. The relation \ref{h01n} implies 
 that $\hat{H}^{(1)}(F, \apr{\rg{L}}) \cong R$ as a cyclic $\F_p[ T ]$ - module. Moreover, applies Tchebotarew to deduce 
 that there is a norm coherent sequence $v \in \apr{\rg{L}} \setminus (p , T) \cdot
  \apr{\rg{L}}$, such that the image $\beta(v) \in \hat{H}^{(1)}(F, \apr{\rg{L}}) $ generates this cyclic $R$-module. Since
  $\hat{H}^{(1)}(F, \apr{\rg{L}})$ is annihilated by $p$ and
  $\nu_{m',1}$, if follows in particular that $T^{d-1} v \not \in p
  \cdot \apr{\rg{L}}$, with $d = \deg(\nu_{m',1})$. We conclude from
  the choice of $m'$, that for $\alpha = v$ and $\K_{aux} = \rg{L}$
  the condition \rf{palpha} is verified.
\end{proof}

\subsection{The working base field and a split Thaine shift}
In this section we shall construct what has been already designed 
as the working base field $\K$ over which we develop
our proof. Suppose that some CM number field $\K_{start}$ is provided,
in which the Leopoldt Conjecture fails for the prime $p$ and let
$\K_{-1} := \K_{start}^{(n)}[ \zeta_p ]$ be the galois closure of
$\K_{start}$, to which we adjoin the \nth{p} roots of unity. Let
$\kappa > 0$ be determined by $\K_{-1} \cap \B_{\infty} =
\B_{\kappa}$. We let then $m = \kappa$ and $k = 3$ in the Lemma
\ref{paux} and let $\K_{aux}$ be the number field, the existence of
which is proved in the Lemma. We then define $\K_{ini} = \K_{-1} \cdot \K_{aux}$ and
by construction, we also have $\K_{ini} \cap \B_{\infty} =
\B_{\kappa}$. Thus, if $\tilde{\gamma} \in \Gal(\K_{\infty}/\Q)$ is a
lift of a topological generator of $\Gal(\K_{aux, \infty}/\K_{aux})$,
then $\tilde{\tau} := \tilde{\gamma}^{p^{\kappa-1}}$ is a topological
generator of $\Gal(\K_{\infty}/\K)$.  Letting $\tilde{T} =
\tilde{\tau}-1$, we see that the Lemma \ref{paux} provides a sequence
$\alpha \in \apr{\K_{aux}}$ such that $\tilde{T}^k \alpha \not \in p
\cdot \apr{\K_{aux}}$.

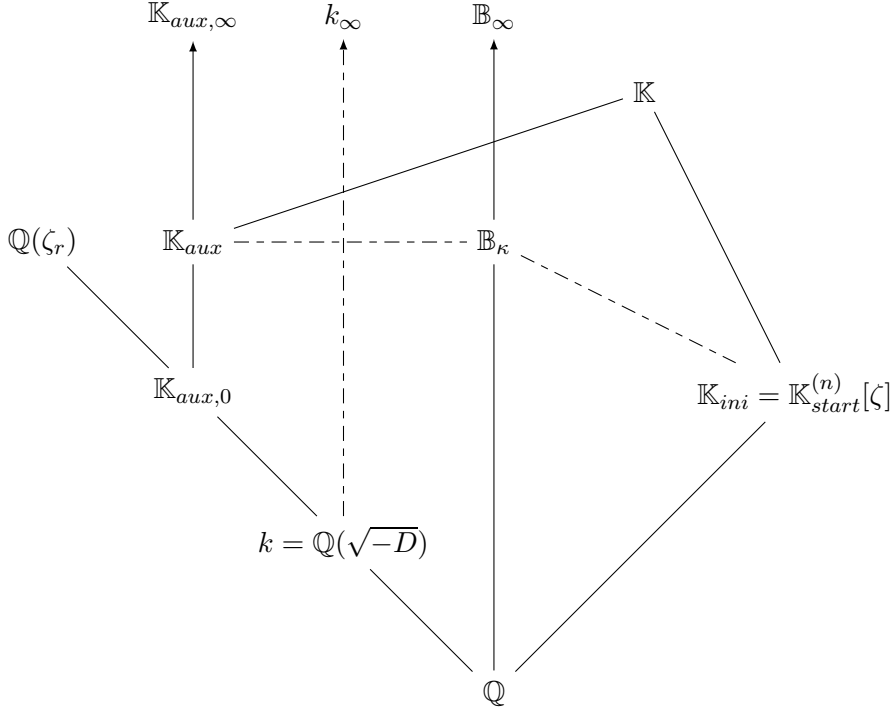
\begin{figure}
\centering
\begin{tikzpicture}[node distance = 2cm, auto]
      \node (Q) {$\mathbb{Q}$};
      \node (k) [above of=Q, left of=Q,node distance = 2cm] {$k=\mathbb{Q}(\sqrt{-D})$};
      \node (Kaux0) [above of=k, left of=k] {$\mathbb{K}_{aux,0}$};
      \node (Cycl) [above of=Kaux0, left of=Kaux0, node distance = 2cm] {$\mathbb{Q}(\zeta_r)$};
      \node (Kaux) [above of=Kaux0, node distance = 2cm] {$\mathbb{K}_{aux}$};
      \node (Bk) [above of=Q, node distance = 6cm] {$\mathbb{B}_{\kappa}$};
      \node (B) [above of=Bk, node distance = 3cm] {$\mathbb{B}_{\infty}$};
      \node (Kstart) [above of=Q, above of=Q, right of=Q, right of=Q, node distance = 2cm] {$\mathbb {K}_{ini} = \mathbb {K}_{start}^{(n)}[ \zeta ]$};
      \node (K) [above of=Kstart, above of=Kstart, left of=Kstart] {$\mathbb{K}$};
      \node (Kauxinfty) [above of=Kaux,  node distance = 3cm] {$\mathbb{K}_{aux,\infty}$};
      \node (kinf) [above of=k, node distance = 7cm] {$k_{\infty}$};
      \draw[-] (Q) to node {} (k);
      \draw[-] (Q) to node {} (Kstart);
      \draw[-] (k) to node {} (Kaux0);
      \draw[-] (Kaux0) to node {} (Cycl);
      \draw[-] (Kaux0) to node {} (Kaux);
      \draw[->,>=latex] (Kaux) to node {} (Kauxinfty);
      \draw[-] (Q) to node {} (Bk);
      \draw[->,>=latex,dash pattern = on 2 mm off 1 mm on 1 mm off 1 mm] (k) to node {} (kinf);      
      \draw[-] (Kaux) to node {} (K);
      \draw[dash pattern = on 2 mm off 1 mm on 1 mm off 1 mm] (Kaux) to node {} (Bk);
      \draw[->,>=latex] (Bk) to node {} (B); 
      \draw[dash pattern = on 2 mm off 1 mm on 1 mm off 1 mm] (Bk) to node {} (Kstart);
      \draw[-] (Kstart) to node {} (K);     
\end{tikzpicture}
\caption{The working base field $\K$} \label{fig:generalconstr}
\end{figure}

Let now $\K := \K_{ini, n_0}$, with $n_0$ determined by the conditions
N1, N2 with respect to the field $\K_{ini}$, let $\kappa =
\kappa(\K)$. Let $N > \kappa$ be such that
$\prk(A^-(\K_N)) = \prk(\apr{\K}) < p^{N-\kappa-1}$.  Consider a prime
ideal $\eu{q} \subset \K_N$ which is totally split and such that
$\Norm_{\K_N/\K_{aux,N}}(\eu{q}) \in \alpha_N$. Such a prime can be
determined by Tchebotarew.  We let $\tilde{a}_N = [ \eu{q}^{1-\jmath}
]$ and $\tilde{a}_m = N_{N,m}(a_{N})$ for $m < N$. We can extend the
sequence $\tilde{a}_m$ to a norm coherent $\tilde{a} \in
\apr{\K_{\infty}}$ with $\Norm_{\K/\K_{aux}}(\tilde{a}) = \alpha$;
this can also be proved using the Tchebotarew theorem.

Let $\M \supset \K_{\infty}$ be the compositum of all CM
$\Z_p$-extensions of $\K$. Since Leopoldt's Conjecture was assumed to
fail, we have $\zprk(\Gal(\M/\K)) \geq 2$, and in fact larger lower
bounds can be gained from the theory of representations. Let $\eu{q}_0
= \eu{q} \cap \K$ be the prime below $\eu{q}$ and let $D(\eu{q}_0)
\subset \Gal(\M/\K)$ be its decomposition group. Since $\eu{q}_0$ does
not ramify, we have $D(\eu{q}_0) = \Z_p$, so there is at least one CM
$\Z_p$-extension of $\K$ in which $\eu{q}_0$ is totally split and
which contains $\K_N$. Let thus $\KL/\K$ be a CM $\Z_p$-extension with
$\KL \cap \K_{\infty} = \K_N$, with $[ \KL_{n} : \KL_N ] = p^{n-N}$
for all $n > N$, in the natural counting induced by setting $\KL_N =
\K_N $. By Lemma \ref{mubd}, we know that $\mu^-(\KL) = 0$.  We let
$\Gamma = \Gal(\KL/\K)$ and $\tau \in \Gamma$ be a topological
generator such that $\tau \vert_{\K_N} = \tilde{\tau}
\vert_{\K_N}$. We also let $T = \tau - 1$ and note that $T$ and
$\tilde{T}$ act identically on groups related to the field $\K_N$ or
its subfields, such as $A^-(\K_n), n \leq N$. Let $\eu{q}_n \subset
\KL_n$ for $n > N$ build a norm coherent sequence of (totally split)
primes above $\eu{q}$. We define $a_n = [ \eu{q}_n^{1-\jmath} ]$ to be
the respective classes and $a = (a_n)_{n \in \N} \in \apr{\KL}$. We
claim that
\begin{eqnarray}
\label{adeco}
a & = & a_{\lambda} \in  \id{L}(\apr{\KL}) \quad \hbox{ and } \\
T^3 a & \not \in & p \cdot \apr{\KL} \quad \hbox{and} \quad T^3 a_{\lambda, n} \not \in p A^-(\K_n) \quad \forall n \leq M. 
\nonumber 
\end{eqnarray}
Indeed, we have $a_N \not \in (p, T) A^-(\K_N) = A^-(\KL_N)$, so a
fortiori $a \not \in (p, T) \apr{\KL}$.  Moreover, Lemma \ref{mubd}
implies that $a \in \id{L}(\apr{\KL}) = \apr{\KL}$.  The second
condition in \rf{adeco} follows from the construction of $\tilde{a}$.

%

Let now $q \in \N$ be the rational prime below $\eu{q}$ and let $\F
\subset \Q[ \zeta_q ]$ be the subfield of absolute degree $p$ in the
\nth{q} cyclotomic extension. We let the \textit{split Thaine shift}
of $\K$ be $\K' = \K \cdot \F$ with its $\Gamma$-extension $\KL' = \KL
\cdot \F$; obviously, $\Gal(\KL'/\K') = \Gal(\KL/\K) = \Gamma$ and $F
:= \Gal(\F/\Q)$ commutes with $\Delta_0 = \Gal(\K/\Q)$. We let $\nu
\in F$ be a generator of this cyclic group of order $p$, write $s =
\nu - 1$ and $\id{N} = \sum_{i=0}^{p-1} \nu^i $. We note that
\[ \id{N} = p + s f(s) = p u(s) + s^{p-1}, \quad \hbox{for some} \quad f \in \Z_p[ s ], u \in (\Z_p[ s ])^{\times} .\]

The primes above $q$ are ramified in $\KL'_n/\KL_n$ for all $n > 0$
and we let $\eu{Q}_n$ be the ramified primes above $\eu{q}_n$ for $n >
0$.  Let $b_n = [ \eu{Q}^{1-\jmath} ]$ be the induced classes and $b =
(b_n)_{n \in \N} \in \apr{\KL'}$ be the corresponding norm coherent
sequence. Ramification implies that $p b_n = a_n$ for all $n > 0$ and
$s b_n = 0$. We call the sequence $b \in \apr{\KL'}$ a {\em Thaine
  lift of $a$}.  Since $p b = a$, it follows that the minimal
annihilator polynomial of $p b$ divides $f(\tilde{T})$ and thus $(p,
f(\tilde{T})) b \in \id{D}(\apr{\KL'})$ is decomposed. Condition
N2 implies that in fact $T b$ is decomposed.
\begin{figure}
\centering
\begin{tikzpicture}[node distance = 1.5cm, auto]
      \node (Q) {$\mathbb{Q}$};
      \node (BK) [above of=Q, node distance = 2cm] {$\mathbb{B}_{\kappa}$};
      \node (B) [above of=BK, node distance = 7cm] {$\mathbb{B}_{\infty}$};
      \node (K) [above of=Q, right of=Q, right of=Q, node distance = 1.5cm] {$\mathbb{K}$};
      \node (F) [above of=K, left of=K, node distance = 1cm] {$\mathbb{F}$};
      \node (Qz) [above of=F, left of=F, node distance = 1cm] {$\mathbb{Q}(\zeta_q)$};
      \node (KN) [above of=K] {$\mathbb{K}_N = \mathbb{L}_N$};
      \node (LprimeN) [above of=F] {$\mathbb{L}^{'}_N$};
      \node (Kinf) [above of=KN, node distance = 6cm] {$\mathbb{K}_{\infty}$};
      \node (Ln) [above of=KN, right of=KN, node distance = 3cm] {$\mathbb{L}_n$};
      \node (Lprimen) [above of=LprimeN, right of=LprimeN, node distance = 3cm] {$\mathbb{L}^{'}_n=\mathbb{F}.\mathbb{L}_n$};
      \node (Linf) [above of=Ln, right of=Ln,node distance = 3cm] {$\mathbb{L}_{\infty}$};
      \node (Lprimeinf) [above of=Lprimen, right of=Lprimen, node distance = 2.5cm] {$\mathbb{L}^{'}_{\infty}=\mathbb{F}.\mathbb{L}_{\infty}$};
      \draw[-] (Q) to node {} (BK);
      \draw[->,>=latex] (BK) to node {} (B);
      \draw[dash pattern = on 2 mm off 1 mm on 1 mm off 1 mm] (BK) to node {} (K);
      \draw[-] (Q) to node {} (K);
      \draw[-] (K) to node {} (F);
      \draw[->,>=latex] (F) to node {} (Qz);
      \draw[-] (F) to node {} (LprimeN);
      \draw[double distance=3pt, line width=1pt] (K) to node {} (KN);
      \draw[->,>=latex] (KN) to node {} (Kinf);
      \draw[-] (LprimeN) to node {} (Lprimen);      
      \draw[-] (LprimeN) to node {} (KN);      
      \draw[-] (KN) to node {} (Ln);      
      \draw[-] (Lprimen) to node {} (Ln);      
      \draw[->,>=latex] (Lprimen) to node {} (Lprimeinf);
      \draw[->,>=latex] (Ln) to node {} (Linf);
\end{tikzpicture}
\caption{The Thaine split shift of $\KL$} \label{fig:thaine}
\end{figure}

\section{Proof of the main Theorem}
The proof of the main Theorem will follow by simple algebraic arguments
based on the previous results and construction, as soon as we establish
the fact that the Thaine lift must be indecomposed under the given 
assumptions.
\subsection{Cohomologic results concerning the Thaine lift}
We keep the notation of the previous chapter and denote the
Tate-cohomologies by
 \[\hat{H}^i(F, X), i = 0, 1, \]
 with $X$ some finite or infinite abelian group on which $F = < \nu >$
 acts.  These are governed by the Hasse Norm Principle and similar
 properties which are ingredients of the proof of Chevalley's Theorem,
 also called the {\em ambig class formula} \cite{La}, Chapter 13, Lemma 4.1.

 Recall that $\eu{q} \in \K_N$ is totally split in $\KL/\Q$ and we
 denoted by $\eu{q}_0$ the prime of $\K$ below $\eu{q}$. The following
 proposition determines the Tate cohomology groups $\hat{H}^i(F,
 \apr{\KL}), i = 0, 1$ and will be proved in the Appendix. We define
 in addition to $\Delta_0$, the galois groups $\Delta =
 \Gal(\K/\rg{K})$ and $\Delta_n = \Gal(\KL_n/\rg{K})$
\begin{proposition}
\label{genHasse}
Notations being like above, let $\Delta = \Gal(\K/\rg{K})$ and $D = | \Delta |$.
For every $n > N$ we have 
\begin{eqnarray}
\label{h01n} 
 \wh{H}^i(F,A^-(\KL_n)) \cong \Lambda[ \Delta ][ T ]/(p, \omega_n), \quad i = 0, 1. 
\end{eqnarray}
In the limit, $\wh{H}^i(F, \apr{\KL})$ are free $\F_p[[ T ]]$-modules
of rank $D$. Letting 
\begin{eqnarray}
\label{betaproj} 
\beta_i : \apr{\KL} \ra \wh{H}^i(F, \apr{\KL})
\end{eqnarray}
 be the natural projections. With this notation, $\beta_0(b)$ generates 
 $\wh{H}^0(F, \apr{\KL})$ as a free $\F_p[[ T ]] [ \Delta ]$-module.  
 In particular, $\Lambda b$ is indecomposed.
\end{proposition}

\subsection{The final argument}
The assumption that the Leopoldt Conjecture fails for some CM field
$\K_{start}$ led us to the construction of a CM $\Z_p$-extension $\KL$
and a Thaine lift $\KL'$ thereof, together with a norm coherent
sequence $a = \apr{\KL}$. This sequence consists of classes determined
by a sequence of totally split primes, and a Thaine lift $b$ thereof
is determined by ramified primes above the ones defining $a_n$.
Moreover, the Hasse Norm Principle shows that, if the construction is
valid -- and thus Leopoldt's Conjecture is false -- then $b$ must be
indecomposed.

We now investigate the reciprocal relations of $a$ and $b$, and,
using decomposition results, will derive a contradiction. 
Since $a = p b \in \id{L}(\apr{\KL})$ is decomposed, the choice of 
$N$ implies
that $T b \in \id{D}(\apr{\KL'})$, so let
\[ T b = b_{\lambda} + b_{\mu}.\] Then $s b_{\lambda} + s b_{\mu} =
T s b = 0$, so
\[ s b_{\lambda}, s b_{\mu} \in \id{L}(\apr{\KL'}) \cap
\id{M}(\apr{\KL'}) = C^-(\KL'). \] Since the finite torsion module
$C^-$ is annihilated by $T^*$, it follows that
\[ T^* b_{\lambda} \in \Ker( s : \id{L}(\apr{\KL'}) \ra
\id{L}(\apr{\KL'})) . \]

Moreover, since $\hat{H}^0(F, \apr{\KL'})$ is a torsion-free $\F_p[[ T
]]$-module, while the image of $b_{\lambda}$ in this module has a
finite annihilator, it follows that $T^* b_{\lambda} \in
\iota_{\KL'/\KL}(\apr{\KL})$. Let thus $T^* b_{\lambda} = \iota(z)$
for some $z \in \apr{\KL}$. We have the following sharpening: let $s
b_{\lambda} = -s b_{\mu} = x \in C^-(\KL')$. Then $\id{N}(x) = 0$, so
$x \in \Ker(s : C^- \ra C^-) = C^-[ p ]$, so in fact $p x = 0$.

On the other hand,
\[ \id{N}(T b) = T a = \left( p \left(1 + \frac{p-1}{2} s \right) +
  O(s^2) \right) \cdot (b_{\lambda} + b_{\mu}).\] Since $s C^-(\KL') =
0$, we also have $s^2 b_{\lambda} = s^2 b_{\mu} = 0$; we have 
seen that $p s b_{\lambda} = p s b_{\mu} = 0$. Thus
$\id{N}(b_{\lambda}) = p b_{\lambda}$ and $\id{N}(b_{\mu}) = p b_{\mu}
\in \id{M}(\apr{\KL})$, so $p b_{\mu} = 0$.  By comparing parts in the
previous identities, we deduce that
\[ T T^* a = \id{N}(T^* b_{\lambda}) = \id{N}(\iota(z)) = p z \quad
\hbox{ and } \quad p b_{\mu} = 0.\] In particular, since $T^* \equiv T
\bmod p$, we deduce that $T^2 a \in p \cdot \apr{\KL}$.  
At level $N$, this implies
\[ a_{ N} \in p \cdot A^-(\K_N), \]
in contradiction with the relation \rf{adeco}. 
This contradiction shows that the above construction is
inconsistent, so the Leopoldt Conjecture must hold for all CM fields
and all odd primes $p$.

\section{Appendix :Proof of Proposition \ref{fantom}}
Let the module $N = \apr{\K_{\infty}}(T^*)$ be defined in the
cyclotomic $\Z_p$-extension of $\K$, and suppose that $\K$ is a
CM-extension with positive Leopoldt defect and containing the \nth{p}
roots of unity. We have mentioned that $\zprk(\Gal(\M^+/\K_{\infty}))
= \id{D}_L(\K)$, a fact which is proved in all text-books. Let $\Omega
\supset \K_{\infty}$ be the maximal $p$-abelian $p$-ramified extension
of $\K_{\infty}$ and
\[ \Omega_E = \bigcup_n \K_n[ E(\K_n)^{1/p^n} ] \subset \Omega. \]

Let $\M \subset \Omega(\K)$ be the compositum of all $\Z_p$-extensions of
$\K$ and let $\KP = \M^- \cap (\KH^- \cap \Omega_E)$. One can build an
explicit map $\rho : E(\K) \otimes_{\Z} \Z_p \ra \overline{E}$ such
that $\Ker(\rho) \sim \Rad(\KP/\K_{\infty})$ and thus
$\zprk(\Gal(\KP/\K_{\infty})) = \id{D}_L(\K)$ while reflection yields
$T^* \Gal(\KP/\K_{\infty}) = 0$.  The extension $\KP$ was for instance
investigated by Jaulent in \cite{Ja}; we denote it \textit{the phantom
  field} associated to the Leopoldt conjecture.

The $p$-ramified, $p$-abelian, real extensions of $\K_{\infty}$ are
obtained as Kummer extensions by taking roots of classes in $A^-$,
according to the point 2. in Remark \ref{shift}. In fact, if $\Omega^+
= \Omega^{Y^-}$ is the fixed field of the minus part of $Y :=
\Gal(\Omega/\K_{\infty})$, we have $\Omega^+ = \overline{\K}_{\infty}[
(A^-)\pinf ]$. Here $\overline{\K}$ indicates that we might have to
adjoin first the roots of some expressions of the type
$\wp/\overline{\wp}$, with $\wp$ a principal prime of $\K_{\infty}$
above $p$.

The galois properties of the Kummer pairing imply more precisely that
$\M^+ \subset \K_{\infty}[ N\pinf ]$, and since $T \Gal(\M^+ /
\K_{\infty} ) = 0$, it follows that $T^* \Rad(\M^+ /\K_{\infty} ) =
0$. Considering $Z := \Gal(\K_{\infty}[ N\pinf ]/\K_{\infty}) \cong
N^{\bullet}$, it follows in fact that $\M^+ = (\K_{\infty}[ N\pinf
])^{T Z}$. It follows by duality that $\Rad(\M^+ /\K_{\infty} ) \cong
N/(T^* N)$. There is an exact sequence of pseudoisomorsphisms:
\[ 1 \ra N[ T^* ] \ra N \ra N \ra N/(T^* N) \ra 1, \]
in which the central map is $T^* : N \ra N$. From this, we deduce 
\begin{eqnarray*}
  \id{D}_L(\K) & = & \zprk(\Gal(\M^+/\K_{\infty})) = \zprk(\Rad(\M^+/\K_{\infty})) \\ & = & \eprk(N/T^* N) =
  \eprk(N[ T^* ]) . 
\end{eqnarray*}
We have thus shown that $\eprk(A^-[ T^* ]) = \id{D}_L(\K)$ and for
each $\KL \subset \M^+$ there is a sequence $a \in A^-(T^*) \setminus
T^* A^-(T^*)$ with $\LK = \K_{\infty}[ a\pinf ]$. This completes the
proof of the Proposition \ref{fantom}.

\section{Proof of Proposition \ref{genHasse}}
Let $n > N$ be fixed and $x \in A^-(\KL'_n)$ be such that $\beta_0(x)
\neq 0$, with $\beta_0$ the projection in \rf{betaproj}.  If $\eu{R}
\in x$, then $\eu{R}^s = (\rho), \rho \in \KL'$, so $\id{N}(\rho) =
(1)$.  Let $\xi = \zeta_{p^N}$ be a primitive \nth{p^N} root of unity,
so $\xi$ generates the $p$ part of the roots of unity $W(\KL_n)$. One
verifies by local class field theory that $\xi \not \in
\id{N}({\KL_n'}^{\times})$.  It follows that, after possibly modifying
$\rho$ by some non primitive root of unity, we may assume that
$\id{N}(\rho^{1-\jmath}) = 1$ and thus $\rho^{1-\jmath} = w^s$, for
some $w \in {\KL'}^{\times}$. Then $\eu{X} := \eu{R}^{1-\jmath}/(w)$
verifies $\eu{X}^s = (1)$ and $[\eu{X} ] = x^2$. The assumption
$\beta_0(x) \neq 0$ implies that $\eu{X}$ is not a lift from $\K$, so
it must be a ramified ideal. Let $B_n = \Lambda[ \Delta ] b_n$ be the
module generated by all classes which contain ramified ideals. We have
thus shown that $\beta_0 : B_n \ra \wh{H}^0(F, A^-(\KL'_n))$ is a
surjective map. In particular $\wh{H}^0(F, A^-(\KL'_n))$ is an
$\F_p$-module of rank $\prk(\wh{H}^0(F, A^-(\KL'_n))) \leq
\prk(B_n)$. Since $\id{N} \vert_{\iota(A^-(\KL_n))} = \cdot p$, the
group $\wh{H}^1(F, A^-(\KL'_n))$ is also an $\F_p$-module. The class
group $A^-(\KL_n)$ is finite, so the Herbrand quotient vanishes and we
have $| \wh{H}^1(F, A^-(\KL_n)) | = | \wh{H}^0(F, A^-(\KL_n)) |$ and
thus
 \begin{eqnarray}
\label{rks} 
\prk( \wh{H}^1 (F, A^-(\KL'_n)) ) = \prk( \wh{H}^0(F, A^-(\KL'_n)) ) \leq \prk(B_n).
\end{eqnarray}

Let $T_n$ denote the image of $T$ in the ring $R_n := \Lambda/(p,
\omega_n)$. We shall prove that
\begin{eqnarray}
\label{mainiso}
\wh{H}^1(F,A^-(\KL'_n)) \cong \F_p[ \Delta ][ T_n ].
\end{eqnarray}
All the claims of the Proposition follow from this fact. Indeed, the
isomorphism $\beta_1(A^-(\KL'_n))) \cong \F_p[ \Delta ][ T_n ]$ is
consequence of \rf{mainiso}. Since the $p$-ranks of $\wh{H}^i(F,
\A^-(\KL'_n))$ are equal by \rf{rks} and $\prk(\F_p[ \Delta ][ T_n ])$
is the number of distinct ramified primes in $q$, thus the upper bound
established for $\prk(\wh{H}^0(F, A^-(\KL'_n))$, it follows that this
rank is maximal. In particular, the ramified primes in $\id{B}_n = \{
g \eu{Q}_n \ : \ g \in \Gal(\KL_n/\rg{K}) \}$ generate classes with
independent images $\beta_0(g b_n)$ and $B_n/p B_n \cong \wh{H}^0(F,
\A^-(\KL'_n)) \cong \F_p[ \Delta ][ T_n ]$ too.  We have
$\varprojlim_n \F_p[ \Delta ][ T_n ] = \F_p[ \Delta ][[ T ]]$, so in
the projective limit, $\wh{H}^i(F, \apr{\KL'})$ are free $\F_p[[ T
]]$-modules of rank $D$.

It is simple now to show that $b$ must be indecomposed. Indeed, since
we established that $\prk(\Lambda b) = \infty$, we certainly have $b
\not \in \id{L}(\apr{\KL'})$. Likewise, if $b \in \id{M}(\apr{\KL'})$,
the relation $a = p b$ would imply that $\prk(\Lambda a_N) =
\deg(\omega_N)$; however $\prk(\Lambda a_N) \leq \lambda^-(\K)$ and we
may choose $N$ large enough to obtain a contradiction. It remains that
$b$ is indecomposed.

The rest of this section is dedicated to the proof of \rf{mainiso}.
Since $\eu{q}_n$ is totally split above $\Q$, it follows that the
residue field is $\rg{K}_n = \id{O}(\KL_n)/\eu{q}_n \cong \F_q$; let
$\rg{L}_n' = \KL'_{\eu{Q}_n}$ and $\rg{L}_n = \KL_{\eu{q}_n} \cong
\Q_q$. Then $\rg{L}_n' \cong \rg{F}$, with $\rg{F}$ the subfield of
degree $p$ of $\Q_q[ \zeta_q ]$. Let $r' \in \F_q^{\times}$ be a
generator and $r = {r'}^{(q-1)/p^N}$ generate the maximal $p$-cycle
contained in $\F_q^{\times}$. Since $\F/\Q_q$ is tamely ramified, it
follows by class field that
\begin{eqnarray}
\label{ndef}
r \not \in \Norm(\rg{F}^{\times}) \quad \hbox{ and } \\ \nonumber
\Gal(\rg{F}/\Q_q) \cong \lan r \ran / \lan r^p \ran. 
\end{eqnarray} 

Let now $P(\M)$ denote the principal ideals of some number field
$\M$. We deduce from the previous consideration the following global
result:
\begin{lemma}
\label{princid}
Notations being like above,
\begin{eqnarray}
\label{h1gen}
\wh{H}^{(1)}( F, A^-(\KL'_n) ) \cong 
P^-(\KL_n)/\id{N}(P^-(\KL'_n)),
\end{eqnarray}
the isomorphism being one of cyclic 
$\F_p[ \Delta_n ]$-modules. 
\end{lemma}
\begin{proof}
  Note that both modules in \rf{h1gen} are annihilated by $p$. In the
  case of $P^-(\KL_n)/\id{N}P^-(\KL'_n)$ this is a direct consequence
  of 
  \[ (P^-(\KL_n))^p = \id{N}(P^-(\KL_n)) \subset \id{N}(P^-(\KL'_n)) .\] 
  Let $\pi_n : P^-(\KL_n) \ra
  P^-(\KL_n)/\id{N}(P^-(\KL'_n))$ denote the natural projection and
  let $a \in \Ker( \id{N} : A^-(\KL'_n) \ra A^-(\KL'_n) )$. Then $p
  u(s) a = - s^{p-1} a$ and thus $p a = - s^{p-1} u^{-1}(s) a$; a
  fortiori $\beta_1(p a) = 0$ for all $a \in \A^-(\KL'_n)$, so $p
  \wh{H}^{(1)}( F, A^-(\KL) ) = 0$. Let now $a \in \Ker(\id{N} :
  A^-(\KL'_n) \ra A^-(\KL'_n))$ and $\eu{A} \in a$ be some prime . Let
  $(\alpha) = \id{N}(\eu{A})$. The principal ideal $\eu{a} :=
  (\alpha/\overline{\alpha}) \in P^-(\KL_n)$ has image $\pi_n(\eu{a})$
  which depends on $a$ but not on the choice of $\eu{A} \in a$. This
  is easily seen by choosing a different ideal $\eu{B} = (x) \eu{A}
  \in a$: then $\id{N}(\eu{B}^{1-\jmath}) = \eu{a} \cdot
  \id{N}(x/\overline{x}) \in \eu{a} \cdot \id{N}(P^-(\KL'_n)$, and
  $\pi_n(\id{N}(\eu{B}^{1-\jmath})) = \pi_n(\id{N}(\eu{A}^{1-\jmath}))
  = \pi_n(\eu{a})$ depends only on $a$. Suppose now that $\eu{a} \in
  \id{N}(P^-(\KL'_n)$. Then there is some $y \in {\KL'}_n^{\times}$
  such that $\id{N}(\eu{A}^{1-\jmath}) = (\id{N}(y))^{1-\jmath}$ and
  thus $\id{N}(\eu{A}/(y))^{1-\jmath} = (1)$. Since $\wh{H}^{(1)}$
  vanishes for ideals, it follows that there is a further ideal
  $\eu{X} \subset \KL'_n$ such that
  \[ \eu{A}^{1-\jmath} = \left((y) \eu{X}^s \right)^{1-\jmath} , \]
  and thus $a^2 \in (A^-(\KL'_n))^s$. But then $\beta_1(a) = 0$. We
  have shown that there is a map $\lambda : \wh{H}^{(1)}( F,
  A^-(\KL'_n) ) \ra P^-(\KL_n)/\id{N}(P^-(\KL'_n))$ defined by the
  sequence of associations $\beta_1(a) \mapsto \eu{a} \mapsto
  \pi_n(\eu{a})$, which is a well defined injective homomorphism of
  $\F_p$-modules.

  In order to show that $\lambda$ is an isomorphism, let $\eu{x} :=
  (x/\overline{x}) \in P^-(\KL_n) \setminus \id{N}(P^-(\KL'_n))$ be a
  principal ideal that is not a norm from $\KL'_n$. Let the Artin
  symbol of $x$ be $\sigma = \lchooses{\KL'_n/\KL_n}{x} \in
  \Gal(\KL'_n/\KL_n)$; since $\KL'_n/\KL_n$ is real, the complex
  conjugation commutes with $\nu$ and we have
  \[ \lchooses{\KL'_n/\KL_n}{\overline{x}} =
  \lchooses{\KL'_n/\KL_n}{x^{\jmath} } = \nu^{\jmath} = \nu . \]
  Consequently, $\lchooses{\KL'_n/\KL_n}{(x/\overline{x})} = 1$ -- we
  may thus choose, by Tchebotarew, a principal prime $(\rho) \subset
  \KL_n$ with $\rho \cong x/\overline{x} \bmod q$, and which is split
  in $\KL'_n/\KL_n$. Let $\eu{R} \subset \KL'_n$ be a prime above
  $(\rho)$ and $R = [ \eu{R}^{1-\jmath} ] \in (A^-)(\KL'_n)$ be its
  class.  We claim that $\beta_1(R) \neq 0$. Assume not, so $\eu{R} =
  (y) \eu{Y}^s$, with $y \in \KL'_n$ and $\eu{Y} \subset \KL'_n$; thus
  $\id{N}(\eu{R}^{1-\jmath}) = (\id{N}(y/\overline{y})) \cong \eu{x}
  \bmod \id{N}(P^-(\KL'_n))$. Since $(\id{N}(y/\overline{y})) \in
  \id{N}(P)^-(\KL'_n)$ by definition, it follows that $\eu{x} =
  (x/\overline{x}) \in \id{N}(P)^-(\KL'_n)$, which contradicts the
  choice of $x$. The map is thus $\lambda$ surjective, and in view of
  the previous result, it is an isomorphism of $\F_p[ \Delta_n
  ]$-modules, while $\F_p[ \Delta_n ] + \F_p[ T_n ][ \Delta ]$. 
  The proof will be completed if we show that at least
  one of the modules is isomorphic to $\F_p[ \Delta_n ]$. Consider
\[ R_0 := \id{O}(\KL_n)^-/( q  \id{O}(\KL_n) \cdot (\id{O}(\KL_n))^p)   \cong 
(\KL_n^{\times}/\id{N}({\KL'}_n^{\times}) .\]
By the Chinese remainder theorem, and in view of 
\rf{ndef}, we have $R_0 \cong \F_p[ \Gal( \KL_n/\Q) ]$ and consequently
$R_1 := R_0^- \cong \F_p[ \Delta_n ]$. 
Finally, for showing that 
\[ (\KL_n^{\times})^-/\id{N}(({\KL'}_n^{\times})^-) \cong
P^-(\KL_n)/\id{N}(P^-(\KL'_n)), \] we note first that an ideal
$(\alpha^{1-\jmath}) \in P^-(\KL_n)$ is determine by
$(\alpha^{1-\jmath}) \in (\KL_n^{\times})^-$, up to roots of unity:
therefore
\[ \prk(P^-(\KL_n)/\id{N}(P^-(\KL'_n))) \geq \prk(
(\KL_n^{\times})^-/\id{N}(({\KL'}_n^{\times})^-) ) -1.  \] 
Recall that the field $\rg{K}$ in Lemma \ref{paux} is an imaginary
quadratic field which contains no \nth{p} roots of unity. It follows that
\[ \prk(P^-(\rg{K})/\id{N}(P^-(\rg{K} \cdot \F))) =
\prk(\rk{K}^-/(\rg{K} \cdot \F)^-) = 1, \] and since
$P^-(\KL_n)/\id{N}(P^-(\KL'_n))$ is a cyclic $\F_p[ \Delta_n
]$-module, it follows that in fact
\[ \prk(P^-(\KL_n)/\id{N}(P^-(\KL'_n))) = \prk(
(\KL_n^{\times})^-/\id{N}(({\KL'}_n^{\times})^-) ), \] so
\[ \wh{H}^1(F, A^-(\KL'_n)) \cong P^-(\KL_n)/\id{N}(P^-(\KL'_n)) \cong
\F_p[ \Delta ][ T_n ], \] which completes the proof.
\end{proof}

\subsection{Avoiding the assumption $\mu = 0$}
The fact that $\mu = 0$ for cyclotomic extensions of $\Q$ has been
proved in \cite{MiMu}, so it is a legitimate assumption, which
simplifies the proof of Leopoldt's conjecture for such fields. We
mentioned however that the proof may be completed without assuming
this fact. We keep the notations introduced above and assume that
$\id{M}(\K) \neq 0$; let $p^j$ be its exponent and let $f(T) \in \Z_p[
T ]$ be the minimal annihilator polynomial of
$\id{L}(\apr{\K})$. Moreover, we let
\[ \id{E}(\K) = \id{E}_{\mu} \oplus \id{E}_{\lambda}, \quad
\id{E}_{\lambda} = \bigoplus_j \Lambda/(f_j^{e_j}), \] be the
elementary $\Lambda$-module associated to $\apr{\K_{\infty}}$ and
$F(T) := \prod_j f_j(T)^{e_j}$ be its characteristic polynomial.

We shall prove the following alternative to Proposition \ref{mubd}:
\begin{proposition}
 \label{labd}
For $N$ sufficiently large and $\KL/\K$ a CM $\Z_p$-extension
with $\KL \supset \K_N$, we have $f(T) \id{L}(\apr{\KL}) = 0 \bmod p^j$, 
and $\exp(\id{M}(\apr{\KL}) = p^j$.
\end{proposition}
Assuming the Proposition, we can replace the condition N2. by
$\omega_{n_0} \in (p^{j+1}, F(T)) \Lambda$ and the reader may
verify that this implies that $T b$ is decomposed and the
rest of the proof is identical. 

\begin{proof}
The proof of the Proposition is based on facts and ideas from
the Thesis of S\"{o}ren Kleine \cite{Kl}. We first review
the facts and notations from this thesis, which are relevant
in our context. Then Kleine proves the following
generalization of a result on stabilization of $\Lambda$-modules
of Fukuda, which applied only for the case $\mu = 0$.

\begin{fact}
\label{kstab}
Let $\K_{\infty}/\K$ be some $\Z_p$-extension\footnote{Note that here
  we break the convention of denoting the cyclotomic $\Z_p$-extension
  by $\K_{\infty}$, and allow this to be an arbitrary extension. The
  application will be for the cyclotomic one, though, hence the
  exception.} in which the primes above $p$ are totally ramified. Let
$F(T)$ be the characteristic polynomial of $\apro(\KL)$ and $g(T) \in
\Lambda$ be a prime which is coprime to $F(T)$. If $n$ is such that
\[ \Lambda(A_n/(g A_n)) \cong \Lambda(A_{n-1}/ g A_{n-1} ), \]
then $\Lambda(A_n/(g A_n)) \cong \Lambda(A_{m}/ g A_{m} )$ for all
$m > n$. 
\end{fact}

The fact naturally adapts to minus parts and it is a consequence of
Theorem 3.6 and Lemma 3.10, p. 51, resp. 56. \cite{Kl}. The next step
is to choose $g(T)$ in order to relate Fact \ref{kstab} to the Iwasawa
invariants of \textit{any} $\Z_p$-extension which contains
$\K_n$. Returning to our usual notation, i.e. $\K_{\infty}$ denoting
the cyclotomic $\Z_p$-extension of $\K$, this is particularly easy in
our case, in which $\apr{\K_{\infty}}$ has no finite
$\Lambda$-submodules. Also, in view of Lemma \ref{almin}, any CM
$\Z_p$-extension $\KL$ with $\KL \cap \K = \K_n$ will have the finite
torsion $\id{F}(\apr{\KL}) = C^- \cong \ZM{p^n}$.  Moreover, $C_n^- =
\Lambda r_n$ with $r = (r_n)_{n \in \N} \in \apr{\K_{\infty}}( T^* )$.

Let $g(T) \in \Z_p[ T ]$ be a distinguished polynomial which is
coprime to $F(T)$; such a polynomial can be found, given the fact that
$u(T)$ draws values from an infinite group, while $F(T)$ has finitely
many zeroes in $\overline{\Q}_p$. Using the Definition 3.40 of the
$g$-ranks $\mbox{rk}_{g}$, the Proposition 3.41 \cite{Kl} implies that
\begin{eqnarray*}
s_{\K}(g) := \mbox{rk}_{g}( \id{E}(\apr{\K_{\infty}}) & = & \mbox{rk}_{g}\apr{\K_{\infty}} \\
s_{\KL}(g) := \mbox{rk}_{g}( \id{E}(\apr{\KL}) & = & \mbox{rk}_{g}\apr{\KL} + 1.
\end{eqnarray*}

For $g(T) = T + p u(T), u \in \Lambda^{\times}$, $u(T)$ draws values
from an infinite group, while $F(T)$ has finitely many zeroes in
$\overline{\Q}_p$, so we may choose $g$ coprime to $F$. Assuming $N$
sufficiently large, the Fact \ref{kstab} implies $s_{\K}(g) = \mu(\K)
+ \lambda(\K) = \mu(\KL) + \lambda(\KL) + 1$. Take now $M > \deg(F)$
and let $g = T^M + p u(T)$ be coprime to $F$. In this case we deduce
by the same argument that $M \mu(\K) + \lambda(\K) = M \mu(\KL) +
\lambda(\KL) + 1$. Combining the two relations, we deduce that
$\mu(\KL) = \mu(\K)$ and $\lambda(\KL) = \lambda(\K) - 1$, the defect
in $\lambda$ being generated by the sequence $r \in \apr{\K}_{\infty}$
which remains stable and becomes torsion in $\KL$. Since $\apr{\KL}$
and $\apr{\K_{\infty}}$ agree at level $n$, it follows that the
minimal polynomials are also congruent modulo the subexponent of
$A^-(\K_N)$. In conclusion, all claims of Proposition \ref{labd} can
be fulfilled for sufficiently large $N$, which completes the proof.
\end{proof}

\begin{Acknow}
  \nonumber This is an alternative approach to several previous
  attempts which used $\lambda$-type rather than $\mu$-type sequences;
  with that approach it was only possible to prove the case of the
  Leopoldt conjecture, in which $p$ is totally split in $\K$. This new
  approach based on $\mu$-sequences and split Thaine shifts
  grew from discussions with S\"oren Kleine, related to his
  PhD thesis \cite{Kl} that concerns Greenberg's Null Space Conjecture. 

  I thank Cornelius Greither for his careful reading of preliminary
  drafts of this version. His remarks and the questions he asked
  during a period of almost one year, helped to substantially improve
  the quality of the paper and eliminate several flaws. 

  I would like to express my gratitude to the precious few who
  assisted earlier attempts with discussions and comments: Grzegorz
  Banaszak, John Coates, Ralf Greenberg, Hendrik Lenstra, Florian
  Pop. Many years ago, Nguyen Thong Do decidedly insisted that I should 
give up my work on Iwasawa theory, since ``you have done once an error, so
by incomplete induction, you will always fail''. By all means such a
statement has worked, if necessary, as a positive incitement at times
of low productivity and doubt.
  
Ina Kersten and the colleagues at the Mathematical Institute of
G\"ottingen have provided during many years of evolution, an
atmosphere of understanding and encouragement, which was supportive
for long time research: to them my sincere gratefulness. In
G\"ottingen, the graduate students in winter term of 2015-2016, Pavel
Coupek, Vlad Cri\c{s}an and Katharina M\"uller presented in a seminar
this work, I owe them for their great help in the final redaction.

  Last but not least, this is to Theres and Seraina, who indulged over
  years with an ambig family presence of the researcher.
\end{Acknow}

\bibliographystyle{abbrv} 
\bibliography{leoMu3}
\end{document}